
\input amstex
\mag=\magstep1
\documentstyle{amsppt}

\hcorrection{3truemm}
\vcorrection{-7.2truemm}

\def\L#1{\Cal L[#1]}

\def\IL#1{{\Cal L}^{-1}[#1]}

\def\RE#1{\operatorname{Re}(#1)}
\def\IM#1{\operatorname{Im}(#1)}

\def\res{\operatorname{Res}}

\CenteredTagsOnSplits

\topmatter
\title
The probability distributions of the first
hitting times of Bessel processes
\endtitle
\leftheadtext{Y. HAMANA and H. MATSUMOTO}
\author
Yuji Hamana {\rm and} Hiroyuki Matsumoto
\endauthor
\rightheadtext{HITTING TIMES OF BESSEL PROCESSES}
\affil
Kumamoto University and Yamagata University
\endaffil
\thanks
Partly supported by the Grant-in-Aid for Scientific Research (C)
No.20540121 and 23540183, Japan Society for the Promotion of Science.
\endthanks
\abstract
We consider the first hitting times of the Bessel processes.
We give explicit expressions for the distribution functions
by means of the zeros of the Bessel functions.
The resulting formula is simpler and easier to treat than
the corresponding one which has been already obtained.
\endabstract
\endtopmatter

\document


\head
1. Introduction.
\endhead

\noindent
In this article we consider the first hitting time of 
the Bessel process, which itself is an interesting object
and is one of the important tools to study several problems
in probability theory. By general theory of one-dimensional
diffusion processes, the Laplace transform of the distribution
satisfies an eigenvalue problem for the generator
and it is given by a ratio of the modified Bessel functions.

Except some special cases it is not easy to invert
the Laplace transforms. When the index $\nu$ of the Bessel process
is a half integer $n+1/2$, $n\in\Bbb N$, the Macdonald function
$K_\nu$ is of a simple form. In this case, it turns that $K_{\nu+1}/K_\nu$
is represented by the ratio of polynomials. With the help of
the partial fraction decomposition, Hamana~\cite{9, 10} recently
has inverted the Laplace transform and applied the results
to show the explicit form and the asymptotic behavior
of the expected volume of the Wiener sausage
for the odd dimensional Brownian motion. The method used in \cite{10}
requires some formulae for the zeros of $K_\nu$.

When the starting point of a Bessel process
is the nearer to the boundary $0$ than its arrival point,
the Laplace transform of the first hitting time is
given by a ratio of the modified Bessel function $I_\nu$. 
In this case, Kent~\cite{14} has given an explicit expression 
(see (2.7), (2,8) and (2.9) below) by means of the zeros of $I_\nu$
or the Bessel function $J_\nu$ for the density in his general framework.

The purpose of this paper is to show explicit formula
for the distribution function in the other case,
where the boundary $\infty$ is natural and 
the Laplace transform is written by a ratio of $K_\nu$'s.
For our expressions we need the zeros of $K_\nu$.
In order to prove the results we represent the ratio of
the Macdonald functions by using contour integrals
of functions easier to treat, and invert the Laplace transform.
Recently Byczkowski et al \cite{2,~3} have given similar
but different expressions for the densities of
the first hitting times, and applied the results
to some study on geometric Brownian motions and
hyperbolic Brownian motions. We use the same curve
for the contour integral. But, we show a decomposition
of the Bessel function ratio and use it, which makes
our expression simpler.

Moreover we should mention
that Ismail has considered several types of functions
given as the ratios of modified Bessel functions and
discussed when the functions are completely monotone.
See \cite{11} and the references therein. In particular
he showed that the function on the right hand side of (2.5) below
is completely monotone by expressing it as a Stieltjes transform
of some function. This result gives another expression
for the distribution function of the first hitting time,
while it is more complicated.

Ismail has also shown that the function $K_{\nu+1}/K_\nu$ is 
completely monotone in a similar way. In our context
such a function appears if we give an expression for
the L{\'e}vy measure of the distribution of the first hitting time.
The L{\'e}vy measure will be deduced in the forthcoming paper.
General theory on the infinitely divisibility of the distributions of
the first hitting times of $1$-dimensional diffusion processes 
is given Yamazato~\cite{19}.

This article is organized as follows. We give the main result, Theorem 2.2,
in the next Section 2 and prove it in Section 3.
Section 4 is devoted to the asymptotic behavior
of the tail probability of the first hitting time,
which is obtained as an application of the result.


\head
2. The first hitting time of the Bessel processes.
\endhead

\noindent
For $\nu\in\Bbb R$ the one-dimensional diffusion process
with infinitesimal generator
$$
\Cal G^{(\nu)}=
\frac12\frac{d^2}{dx^2}+\frac{2\nu+1}{2x}\frac d{dx}
=\frac1{2x^{2\nu+1}}\frac d{dx}\biggl(x^{2\nu+1}\frac d{dx}\biggr)
$$
is called the Bessel process with index $\nu$.
If $2\nu+2$ is a positive integer, the Bessel process
is identical in law with the radial motion of a
$(2\nu+2)$-dimensional Brownian motion. Hence, $2\nu+2$ is called
the dimension of the Bessel process.

The classification of boundary points gives the following information.
The endpoint $\infty$ is a natural boundary for any $\nu\in\Bbb R$.
For $\nu\geqq0$, $0$ is an entrance and not exit boundary.
For $-1<\nu<0$, $0$ is a regular boundary, which is instantly
reflecting. For $\nu\leqq-1$, $0$ is an exit but not entrance boundary.
For more details, see \cite{12, 17} for example.

For $a,b\geqq0$ we denote by $\tau_{a,b}^{(\nu)}$ the first hitting time
to $b$ of the Bessel process with index $\nu$ starting at $a$. 
By general theory of one-dimensional diffusion processes,
we can evaluate the Laplace transform of the distribution of
$\tau_{a,b}^{(\nu)}$ by solving an eigenvalue problem.
In fact, the function
$$
x\mapsto E[e^{-\lambda\tau_{x,b}^{(\nu)}}]
$$
is increasing (decreasing) on $[0,b)$ (resp. $(b,\infty)$) and
satisfies
$$
\Cal G^{(\nu)}u=\lambda u,\quad u(b)=1.
$$

The following explicit expressions for $E[e^{-\lambda\tau_{a,b}^{(\nu)}}]$
are well known (cf. \cite{7, 13}): for $\lambda>0$, if $b>0$ and $\nu>-1$,
$$
E[e^{-\lambda\tau_{0,b}^{(\nu)}}]=
\frac{(b\sqrt{2\lambda})^\nu}{2^\nu\varGamma(\nu+1)}
\frac1{I_\nu(b\sqrt{2\lambda})};
\tag2.1
$$
if $0<a\leqq b$ and $\nu>-1$,
$$
E[e^{-\lambda\tau_{a,b}^{(\nu)}}]=
\frac{a^{-\nu}I_\nu(a\sqrt{2\lambda})}
{b^{-\nu}I_\nu(b\sqrt{2\lambda})};
\tag2.2
$$
if $0<a\leqq b$ and $\nu\leqq-1$,
$$
E[e^{-\lambda\tau_{a,b}^{(\nu)}}]=
\frac{a^{-\nu}I_{-\nu}(a\sqrt{2\lambda})}
{b^{-\nu}I_{-\nu}(b\sqrt{2\lambda})};
\tag2.3
$$
if $a>0$ and $\nu<0$,
$$
E[e^{-\lambda\tau_{a,0}^{(\nu)}}]=
\frac{2^{\nu+1}}{\varGamma(|\nu|)(a\sqrt{2\lambda})^\nu}
K_\nu(a\sqrt{2\lambda});
\tag2.4
$$
if $0<b\leqq a$ and $\nu\in\Bbb R$,
$$
E[e^{-\lambda\tau_{a,b}^{(\nu)}}]=
\frac{a^{-\nu}K_\nu(a\sqrt{2\lambda})}
{b^{-\nu}K_\nu(b\sqrt{2\lambda})}.
\tag2.5
$$
Here $\varGamma$ is the gamma function
and $I_\nu$ and $K_\nu$ denote modified Bessel functions
of the first and the second kinds of order $\nu$, respectively.
Both $I_\nu$ and $K_\nu$ are the
solutions of the modified Bessel differential equation
$$
z^2\frac{d^2w}{dz^2}+z\frac{dw}{dz}-(z^2+\nu^2)w=0.
\tag2.6
$$

The distribution function of $\tau_{a,b}^{(\nu)}$ for $a<b$
was obtained in \cite{14}. If $b>0$ and $\nu>-1$,
$$
P(\tau_{0,b}^{(\nu)}\leqq t)=1-\frac1{2^{\nu-1}\varGamma(\nu+1)}
\sum_{k=1}^\infty\frac{j_{\nu,k}^{\nu-1}}{J_{\nu+1}(j_{\nu,k})}
e^{-\frac{j_{\nu,k}^2}{2b^2}t}.\tag2.7
$$
If $0<a<b$, we have that, $\nu>-1$
$$
P(\tau_{a,b}^{(\nu)}\leqq t)=1-2\biggl(\frac ba\biggr)^\nu\sum_{k=1}^\infty
\frac{J_\nu(aj_{\nu,k}/b)}
{j_{\nu,k}J_{\nu+1}(j_{\nu,k})}e^{-\frac{j_{\nu,k}^2}{2b^2}t}
\tag2.8
$$
and that, for $\nu\leqq-1$
$$
P(\tau_{a,b}^{(\nu)}\leqq t)=\biggl(\frac ba\biggr)^{2\nu}
-2\biggl(\frac ba\biggr)^\nu\sum_{k=1}^\infty
\frac{J_{-\nu}(aj_{-\nu,k}/b)}
{j_{-\nu,k}J_{-\nu+1}(j_{-\nu,k})}e^{-\frac{j_{-\nu,k}^2}{2b^2}t}.
\tag2.9
$$
Here $J_\mu$ is the Bessel function of the first kind
of order $\mu$ and $\{j_{\nu,k}\}_{k=1}^\infty$ is
the increasing sequence of positive zeros of $J_\nu$.

When $b>0$ and $2\nu+2$ is a positive integer,
Ciesielski and Taylor~\cite{4} have already shown (2.7).
When $2\nu+2$ is not a positive integer, it is possible to derive
the inverse Laplace transform of the right hand side of (2.1)
by the same methods as those they used to prove
Theorem 1 in \cite{4}.
The formula (2.8) immediately shows that, if $0<a<b$ and $\nu>-1$,
$$
P(\tau_{a,b}^{(\nu)}>t)=2\biggl(\frac ba\biggr)^\nu
\frac{J_\nu(aj_{\nu,1}/b)}{j_{\nu,1}J_{\nu+1}(j_{\nu,1})}
e^{-\frac{j_{\nu,1}^2}{2b^2}t}\{1+o(1)\}.
$$
Similar asymptotic result in the case where $a=0$ and 
$2\nu+2$ is an integer greater than 2 (Brownian case) was used in
\cite{4}
to show the law of iterated logarithm for the total time 
spent by the Bessel process in $(0,b)$ as $b\downarrow0$.

\medskip

\noindent
{\bf Remark 2.1.} We find in \cite{5}, p.104 the formula
$$
\frac{J_\nu(cz)}{J_\nu(z)}=c^\nu+\sum_{k=1}^\infty
\frac{2z^2}{j_{\nu,k}(j_{\nu,k}^2-z^2)}
\frac{J_\nu(c j_{\nu,k})}{J_{\nu+1}(j_{\nu,k})}
$$
for $0\leqq c\leqq 1$, which implies (2.8) and (2.9)
by virtue of $I_\nu(x)=e^{-i\pi\nu/2}J_\nu(xe^{i\pi/2})$
for $x>0$ (cf. \cite{18}, p.77).

\medskip

In the case $a>0$ and $\nu<0$, we can easily check 
$$
P(\tau_{a,0}^{(\nu)}\leqq t)=\frac{2^\nu}{\varGamma(|\nu|)a^2}
\int_0^t s^{\nu-1}e^{-\frac{a^2}{2s}},
$$
by (2.4) and the formula
$$
K_\nu(z)=\frac12\biggl(\frac z2\biggr)^\nu
\int_0^\infty e^{-t-\frac{z^2}{4t}}t^{-\nu-1}dt.
$$

To give our result on the distribution functions of $\tau_{a,b}^{(\nu)}$
in the case of $0<b<a$, we need to recall some facts
about the zeros of the Bessel function $K_\nu$. For $\nu\in\Bbb R$
we denote by $N(\nu)$ the number of zeros of $K_\nu$. It is known
that $N(\nu)=|\nu|-1/2$ if $\nu-1/2$ is an integer
and that $N(\nu)$ is the even number closest to $|\nu|-1/2$ otherwise.
We remark that $N(\nu)=0$ if $|\nu|<3/2$ and $N(\nu)\geqq1$
if $|\nu|\geqq3/2$.
Each zero, if exists, lies in the half plain
$\{z\in{\Bbb C}\,;\,\RE{z}<0\}$, denoted by ${\Bbb C}^-$.
In this case, we write $z_{\nu,1},z_{\nu,2},\dots,z_{\nu,N(\nu)}$
for the zeros. Since $K_\nu$ is a solution of (2.6),
all zeros of $K_\nu$ are of multiplicity one
by the uniqueness of the solution of ordinary differential equations.
This means that all zeros of $K_\nu$ are distinct.
If $\nu-1/2$ is not an integer, there are no real zeros.
For details, see \cite{18}, pp.511--513. 

\medskip

\proclaim{Theorem 2.2}
Let $0<b<a$. For $\mu\geqq0$ and $c>1$, we set
$$
L_{\mu,c}(x)=
\frac{\cos(\pi\mu)\{I_\mu(cx)K_\mu(x)-I_\mu(x)K_\mu(cx)\}}
{\{K_\mu(x)\}^2+\pi^2\{I_\mu(x)\}^2+2\pi\sin(\pi\mu)K_\mu(x)I_\mu(x)}.
$$
\text{\rm{(1)}} If $\nu=\pm1/2$,
$$
P(\tau_{a,b}^{(\nu)}\leqq t)=\biggl(\frac ba\biggr)^{\nu+|\nu|}
\int_0^t \frac{a-b}{\sqrt{2\pi s^3}}
e^{-\frac{(a-b)^2}{2s}}ds.
$$
\text{\rm{(2)}} If $|\nu|<3/2$ and $\nu\neq\pm1/2$,
$$
\split
P(\tau_{a,b}^{(\nu)}\leqq t)=
&\biggl(\frac ba\biggr)^{\nu+|\nu|}
\int_0^t \frac{a-b}{\sqrt{2\pi s^3}}
e^{-\frac{(a-b)^2}{2s}}ds\\
&-\biggl(\frac ba\biggr)^\nu \int_0^t \frac{a-b}{\sqrt{2\pi s^3}}
e^{-\frac{(a-b)^2}{2s}}
\biggl[\int_0^\infty \frac{L_{|\nu|,a/b}(x)}x
e^{-\frac{x(a-b)\sqrt t}{b\sqrt s}}dx\biggr]ds.
\endsplit
$$
\text{\rm{(3)}} If $\nu-1/2$ is an integer and $\nu\neq\pm1/2$,
$$
\split
P(\tau_{a,b}^{(\nu)}\leqq t)=
&\biggl(\frac ba\biggr)^{\nu+|\nu|}
\int_0^t \frac{a-b}{\sqrt{2\pi s^3}}
e^{-\frac{(a-b)^2}{2s}}ds\\
&-\biggl(\frac ba\biggr)^\nu \sum_{j=1}^{N(\nu)}
\frac{K_\nu(az_{\nu,j}/b)}{z_{\nu,j}K_{\nu+1}(z_{\nu,j})}
\int_0^t \frac{a-b}{\sqrt{2\pi s^3}}
e^{-\frac{(a-b)^2}{2s}+\frac{z_{\nu,j}(a-b)\sqrt t}{b\sqrt s}}ds.
\endsplit
$$
\text{\rm{(4)}} If $\nu-1/2$ is not an integer and $|\nu|>3/2$,
$$
\split
P(\tau_{a,b}^{(\nu)}\leqq t)=
&\biggl(\frac ba\biggr)^{\nu+|\nu|}
\int_0^t \frac{a-b}{\sqrt{2\pi s^3}}
e^{-\frac{(a-b)^2}{2s}}ds\\
&-\biggl(\frac ba\biggr)^\nu \sum_{j=1}^{N(\nu)}
\frac{K_\nu(az_{\nu,j}/b)}{z_{\nu,j}K_{\nu+1}(z_{\nu,j})}
\int_0^t \frac{a-b}{\sqrt{2\pi s^3}}
e^{-\frac{(a-b)^2}{2s}+\frac{z_{\nu,j}(a-b)\sqrt t}{b\sqrt s}}ds\\
&-\biggl(\frac ba\biggr)^\nu \int_0^t \frac{a-b}{\sqrt{2\pi s^3}}
e^{-\frac{(a-b)^2}{2s}}
\biggl[\int_0^\infty \frac{L_{|\nu|,a/b}(x)}x
e^{-\frac{x(a-b)\sqrt t}{b\sqrt s}}dx\biggr]ds.
\endsplit
$$
\endproclaim

\medskip

From this theorem we can deduce the asymptotic behavior of
$P(\tau_{a,b}^{(\nu)}>t)$ for $0<b<a$, which will be discussed
in Section 4. It should be mentioned
that, when $\nu=0$, it is obtained from (2.5). In fact, we have
$$
\split
\int_0^\infty e^{-\lambda t}P(\tau_{a,b}^{(\nu)}>t)dt
&=\frac1\lambda\frac{K_0(b\sqrt{2\lambda})-K_0(a\sqrt{2\lambda})}
{K_0(b\sqrt{2\lambda})}\\
&=\frac{2\log(a/b)}{\lambda\log(1/\lambda)}\{1+o(1)\}
\endsplit
\tag2.10
$$
as $\lambda\to0$.
The Tauberian theorem of the exponential type immediately yields
$$
P(\tau_{a,b}^{(\nu)}>t)=\frac{2\log(a/b)}{\log t}
+o\biggl(\frac1{\log t}\biggr)
\tag2.11
$$
as $t\to\infty$ (cf. \cite{6}, p.446).
In order to derive the second equality of (2.10), we have applied
the asymptotic behavior of $K_0(x)$ as $x\downarrow0$ (cf. (3.15)).
We can directly deduce (2.11) from Theorem 2.2
without the Tauberian theorem, however, the calculation is left
to the reader.

In case of $\nu\neq0$,
we can not obtain the convenient formula like (2.10)
which admits us to apply the Tauberian theorem
in a straightforward way.


\head
3. The distribution function of the first hitting time.
\endhead

\noindent
From now on, for a suitable function $f$,
the notation $\L f$ implies the Laplace transform of $f$
and the inverse Laplace transform of $f$ is denoted by $\IL f$.
All formulae concerning Laplace and inverse Laplace transforms
can be found in \cite{16}.

This section is devoted to a proof of Theorem 2.2.
For $t>0$ and $\nu\in\Bbb R$ let
$$
F_{a,b}^{(\nu)}(t)=P(\tau_{a,b}^{(\nu)}\leqq t).
$$
A standard formula shows that
$$
\L{F_{a,b}^{(\nu)}}(\lambda)=
\frac1\lambda E[e^{\lambda\tau_{a,b}^{(\nu)}}]
$$
for $\lambda>0$. For simplicity
we put $G_{a,b}^{(\nu)}(t)=F_{a,b}^{(\nu)}(2b^2 t)$. Then we have
$$
\L{G_{a,b}^{(\nu)}}(\lambda)
=\frac1{2b^2}\L{F_{a,b}^{(\nu)}}\biggl(\frac\lambda{2b^2}\biggr).
$$
Set $\alpha=a/b>1$. It follows from (2.5) that, for $\lambda>0$
$$
\L{G_{a,b}^{(\nu)}}(\lambda)=
\dfrac1{\alpha^\nu}\dfrac{K_\nu(\alpha\sqrt{\lambda})}
{\lambda K_\nu(\sqrt{\lambda})}.
$$
Since $K_\nu=K_{-\nu}$ for $\nu\geqq0$, it is sufficient to consider
the case where $\nu\geqq0$.

Let $\nu\geqq0$ and $c>1$. We first assume that
$\nu-1/2$ is an integer. In this case,
there is a suitable polynomial $\psi_\nu$ of order $\nu-1/2$
on $\Bbb C$ such that $\psi_\nu(0)\neq0$ and
$$
z^\nu K_\nu(z)=\sqrt{\frac\pi2}e^{-z}\psi_\nu(z)
\tag3.1
$$
for $z\in D$ (cf. \cite{15, 18}). Here the notation $D$ has been
used to denote the set of points $z\in\Bbb C\setminus\{0\}$
with $|\arg z|<\pi$. For example,
$$
\psi_{1/2}(z)=1,\quad \psi_{3/2}(z)=1+z,
\quad \psi_{5/2}(z)=3+3z+z^2.
$$
The function $z^\nu K_\nu(z)$ is extended
to a entire function and all zeros of $\psi_\nu$
are the same as those of $K_\nu$. For $z\in\Bbb C$ let
$$
\psi_{\nu,c}(z)=\frac{e^{-(c-1)z}\psi_\nu(cz)}{c^\nu\psi_\nu(z)}.
$$
Then $\psi_{\nu,c}$ is a single-valued meromorphic function on $\Bbb C$
and it holds that
$$
\psi_{\nu,c}(z)=\frac{K_\nu(cz)}{K_\nu(z)}
$$
for $z\in D$. Therefore, if $z\in\Bbb C$ is not a zero of $K_\nu$,
we have
$$
\psi_{\nu,c}(z)=\lim_{v\to z}\frac{K_\nu(cv)}{K_\nu(v)},
$$
which implies that $K_\nu(cx)/K_\nu(x)$
can be determined uniquely for $x<0$
if $x$ is not a zero of $K_\nu$.

Recall our notation $z_{\nu,1},\dots, z_{\nu,N(\nu)}$
for the zeros of $K_\nu$. Let $w$ is a point in $D$ with $K_\nu(w)\neq0$.
We take $R$ so large that $w$ and all zeros of $K_\nu$
are inside $C(R)$, a circle whose center is the origin
and radius R.

We set
$$
\varTheta(R)=\frac1{2\pi i}\int_{C(R)}g_{\nu,c}^w(z)dz,
\tag3.2
$$
where
$$
g_{\nu,c}^w(z)=\frac{we^{(c-1)z}\psi_{\nu,c}(z)}{z(z-w)}
$$
for $z\in\Bbb C$. The singular points of $g_{\nu,c}^w$ are
$0$, $w$ and zeros of $K_\nu$,
which are all poles of order $1$. The residue theorem yields that,
if $N(\nu)=0$,
$$
\varTheta(R)=\res(0;g_{\nu,c}^w)+\res(w;g_{\nu,c}^w)
$$
and that, if $N(\nu)\geqq1$,
$$
\varTheta(R)=\res(0;g_{\nu,c}^w)+\res(w;g_{\nu,c}^w)
+\dsize\sum_{j=1}^{N(\nu)}\res(z_{\nu,j};g_{\nu,c}^w),
$$
where $\res(v;f)$ is the residue of a function $f$ at a pole $v$.
By definition of the function $g_{\nu,c}^w$, we have
$$
\res(0;g_{\nu,c}^w)=-\psi_{\nu,c}(0)=-\frac1{c^\nu}
$$
and
$$
\res(w;g_{\nu,c}^w)=e^{(c-1)w}\psi_{\nu,c}(w)
=e^{(c-1)w}\frac{K_\nu(cw)}{K_\nu(w)}.
$$
If $N(\nu)\geqq1$, the residue of $g_{\nu,c}^w$ at $z_{\nu,j}$
is equal to
$$
\lim_{z\to z_{\nu,j}}\frac{w\psi_\nu(cz)}{c^\nu z(z-w)}
\frac{z-z_{\nu.j}}{\psi_\nu(z)}=
\frac{w\psi_\nu(cz_{\nu,j})}{c^\nu z_{\nu,j}(z_{\nu,j}-w)\psi_\nu'(z_{\nu,j})}
$$
for $1\leqq j\leqq N(\nu)$. Since Lemma 3.1 in \cite{10} gives that,
if $\psi_\nu(z_0)=0$,
$$
\psi_{\nu+1}(z_0)=-z_0\psi_\nu'(z_0),
$$
we obtain
$$
\res(z_{\nu,j};g_{\nu,c}^w)=
-\frac{w\psi_\nu(cz_{\nu,j})}{c^\nu(z_{\nu,j}-w)\psi_{\nu+1}(z_{\nu,j})}.
\tag3.3
$$
If $z\in D$, it follows from (3.1) that
$$
\frac{K_\nu(cz)}{K_{\nu+1}(z)}=
\frac{ze^{-(c-1)z}\psi_\nu(cz)}{c^\nu \psi_{\nu+1}(z)}.
\tag3.4
$$
Then $K_\nu(cz)/K_{\nu+1}(z)$ is extended to a meromorphic
function on $\Bbb C$. This implies that $K_\nu(cx)/K_{\nu+1}(x)$
can be determined uniquely for $x<0$ with $K_{\nu+1}(x)\neq0$.
From (3.3) and (3.4) we deduce
$$
\res(z_{\nu,j};g_{\nu,c}^w)=
-\frac{we^{(c-1)z_{\nu,j}}}{z_{\nu,j}(z_{\nu,j}-w)}
\frac{K_\nu(cz_{\nu,j})}{K_{\nu+1}(z_{\nu,j})}
$$
and
$$
\varTheta(R)=-\frac1{c^\nu}+e^{(c-1)w}\frac{K_\nu(cw)}{K_\nu(w)}
-\sum_{j=1}^{N(\nu)}\frac{we^{(c-1)z_{\nu,j}}}{z_{\nu,j}(z_{\nu,j}-w)}
\frac{K_\nu(cz_{\nu,j})}{K_{\nu+1}(z_{\nu,j})}.
$$
$\varTheta(R)$ tends to $0$ as $R\to\infty$
since $g_{\nu,c}^w(z)=O(|z|^{-2})$.
Hence we obtain
$$
\frac{K_\nu(cw)}{K_\nu(w)}=
\cases
\dfrac{e^{-(c-1)w}}{c^\nu}&\text{if $\nu=\dfrac12$,}\\
\dfrac{e^{-(c-1)w}}{c^\nu}+\dsize\sum_{j=1}^{N(\nu)}
\dfrac{we^{(c-1)(z_{\nu,j}-w)}}{z_{\nu,j}(z_{\nu,j}-w)}
\dfrac{K_\nu(cz_{\nu,j})}{K_{\nu+1}(z_{\nu,j})}\quad
&\text{if $\nu\neq\dfrac12$}
\endcases
\tag3.5
$$
in the case where $\nu-1/2$ is a non-negative integer.

We next consider the case where $\nu-1/2$ is not an integer
and look for a nice expression for $K_\nu(cw)/K_\nu(w)$
like (3.5). If $\nu$ is not an integer, it is well-known
(cf. \cite{18}, p.80) that
$$
K_\nu(ze^{im\pi})=e^{-im\pi\nu}K_\nu(z)
-i\pi\frac{\sin(m\pi\nu)}{\sin(\pi\nu)}I_\nu(z)
$$
for $z\in D$ and $m\in\Bbb Z$. When $\nu$ is an integer,
we also have
$$
K_\nu(ze^{im\pi})=e^{-im\pi\nu}K_\nu(z)-i\pi m(-1)^{(m-1)\nu}I_\nu(z)
$$
for $z\in D$ and $m\in\Bbb Z$, which is easily seen from
$$
\lim_{\mu\to n}K_\mu(z)=K_n(z)
$$
for each integer $n$. Especially, for $z\in D$, we have
$$
\align
&K_\nu(ze^{i\pi})=e^{-i\pi\nu}K_\nu(z)-i\pi I_\nu(z),
\tag3.6\\
&K_\nu(ze^{-i\pi})=e^{i\pi\nu}K_\nu(z)+i\pi I_\nu(z).
\tag3.7
\endalign
$$
It follows from these identities that the function $K_\nu(cz)/K_\nu(z)$
can not be extended to a meromorphic function on $\Bbb C$.
For $z\in D$ let
$$
h_{\nu,c}^w(z)=\frac{we^{(c-1)z}}{z(z-w)}\frac{K_\nu(cz)}{K_\nu(z)}.
$$
In order to give a formula for $K_\nu(cw)/K_\nu(w)$ like (3.5),
we consider the integral of $h_{\nu,c}^w$ on a suitable contour.
However we can not adopt a circle as the contour like (3.2)
since $h_{\nu,c}^w$ can not extend to a meromorphic function on $\Bbb C$.

Let $\varepsilon$ and $R$ be positive numbers with
$2\varepsilon<R$. We set
$$
\theta_{R,\varepsilon}=\operatorname{Arcsin}\frac\varepsilon R.
$$
As a contour, we take the curve $\gamma$ defined by
$$
\align
&\gamma_0\,:\,z=Re^{i\theta},
\,\,-\pi+\theta_{R,\varepsilon}\leqq\theta\leqq\pi-\theta_{R,\varepsilon},\\
&\gamma_1\,:\,z=x+i\varepsilon,
\,\,-R\cos\theta_{R,\varepsilon}\leqq x\leqq0\\
&\gamma_2\,:\,z=\varepsilon e^{i\theta},
\,\,-\pi/2\leqq\theta\leqq\pi/2,\\
&\gamma_3\,:\,z=x-i\varepsilon,
\,\,-R\cos\theta_{R,\varepsilon}\leqq x\leqq0\\
&\gamma=\gamma_0+\gamma_1-\gamma_2-\gamma_3.
\endalign
$$
We take $R$ so large and $\varepsilon$ so small
that $w$ and all zeros of $K_\nu$ are inside $\gamma$.
Then, setting
$$
\varPi(R,\varepsilon)=\frac1{2\pi i}
\int_\gamma h_{\nu,c}^w(z)dz,\quad
\varPi_\ell=\frac1{2\pi i}\int_{\gamma_\ell} h_{\nu,c}^w(z)dz
$$
for $0\leqq \ell\leqq3$, we have
$$
\varPi(R,\varepsilon)=\varPi_0+\varPi_1-\varPi_2-\varPi_3.
$$

The residue theorem yields that, if $N(\nu)=0$,
$$
\varPi(R,\varepsilon)=\res(w;h_{\nu,c}^w)
\tag3.8
$$
and that, if $N(\nu)\geqq1$,
$$
\varPi(R,\varepsilon)=\res(w;h_{\nu,c}^w)+\sum_{j=1}^{N(\nu)}
\res(z_{\nu,j};h_{\nu,c}^w).
\tag3.9
$$
It is obvious that
$$
\res(w;h_{\nu,c}^w)=e^{(c-1)w}\frac{K_\nu(cw)}{K_\nu(w)}.
$$
When $N(\nu)\geqq1$, by using the formula
$$
zK_\nu'(z)-\nu K_\nu(z)=-zK_{\nu+1}(z)
$$
(cf. \cite{18}, p.29), we can show
$$
\res(z_{\nu,j};h_{\nu,c}^w)=
\frac{we^{(c-1)z_{\nu,j}}}{z_{\nu,j}(z_{\nu,j}-w)}
\frac{K_\nu(cz_{\nu,j})}{K_\nu'(z_{\nu,j})}
=-\frac{we^{(c-1)z_{\nu,j}}}{z_{\nu,j}(z_{\nu,j}-w)}
\frac{K_\nu(cz_{\nu,j})}{K_{\nu+1}(z_{\nu,j})}.
$$
It follows  from (3.8) and (3.9) that, if $N(\nu)=0$,
$$
\varPi(R,\varepsilon)=e^{(c-1)w}\dfrac{K_\nu(cw)}{K_\nu(w)}
$$
and that, if $N(\nu)\geqq1$,
$$
\varPi(R,\varepsilon)=
e^{(c-1)w}\dfrac{K_\nu(cw)}{K_\nu(w)}
-\dsize\sum_{j=1}^{N(\nu)}
\frac{we^{(c-1)z_{\nu,j}}}{z_{\nu,j}(z_{\nu,j}-w)}
\frac{K_\nu(cz_{\nu,j})}{K_{\nu+1}(z_{\nu,j})}.
$$

In order to consider asymptotic behavior of
$\varPi(R,\varepsilon)$ as $R\to\infty$,
we need the asymptotic behavior of $K_\nu(z)$ as $|z|\to\infty$.
It is known that, if $|\arg z|<3\pi/2$,
$$
K_\nu(z)=\sqrt{\frac\pi{2z}}e^{-z}\{1+O(|z|^{-1})\}
\tag3.10
$$
(cf. \cite{18}, p.202). However it is not sufficient
since the error term may not be uniform for $\arg z$.
We can improve (3.10) in the following way.

\medskip

\proclaim{Lemma 3.1}
Let $\delta>0$ be given. We have that, for $|\arg z|\leqq 3\pi/2-\delta$
$$
K_\nu(z)=\sqrt{\frac\pi{2z}}e^{-z}\{1+E_1(z)\},
$$
where $|E_1(z)|\leqq C_1/|z|$ and $C_1$ is a positive constant
independent of $z$.
\endproclaim
\demo{Proof}
It is known that, for $-\pi+\delta\leqq \arg w\leqq 2\pi-\delta$,
$$
H_\nu^{(1)}(w)=\sqrt{\frac 2{\pi w}}
e^{i(w-\pi\nu/2-\pi/4)}\{1+E_2(w)\}
$$
holds for the Hankel function $H_\nu^{(1)}$,
where $|E_2(z)|\leqq C_2/|z|$ and $C_2$ is a positive constant
which is independent of $z$. See \cite{15}, p.121 and \cite{18}, p.197.

Let $|\arg z|\leqq3\pi/2-\delta$.
Then we have $-\pi+\delta\leqq \arg(ze^{i\pi/2})\leqq 2\pi-\delta$, and
thus the formula
$$
K_\nu(z)=\frac12 i\pi e^{i\pi\nu/2}H_\nu^{(1)}(ze^{i\pi/2})
$$
(cf. \cite{18}, p.77) immediately yields the assertion of this lemma.
\qed
\enddemo

\medskip

By virtue of Lemma 3.1, we get
$$
|h_{\nu,c}^w(z)|\leqq\frac{|w|}{\sqrt c|z|\cdot |z-w|}
\frac{1+C_1/c|z|}{1-C_1/|z|}
$$
if $z=Re^{i\theta}\in D$ and
$$
|\varPi_0|\leqq\frac1{2\pi}\int_{-\pi+\theta_{R,\varepsilon}}^{\pi-\theta_{R,\varepsilon}}
|h_{\nu,c}^w(Re^{i\theta})|Rd\theta
\leqq\frac{|w|}{\sqrt c(R-|w|)}\frac{1+C_1/cR}{1-C_1/R},
$$
which tends to $0$ as $R\to\infty$.

For the integral $\varPi_1$, we have
$$
\align
\varPi_1
&=\frac w{2\pi i}\int_{-R\cos\theta_{R,\varepsilon}}^0
\frac{e^{(c-1)(x+i\varepsilon)}}{(x+i\varepsilon)(x+i\varepsilon-w)}
\frac{K_\nu(c(x+i\varepsilon))}{K_\nu(x+i\varepsilon)}dx\\
&=\frac w{2\pi i}\int_0^{R\cos\theta_{R,\varepsilon}}
\frac{e^{(c-1)(-x+i\varepsilon)}}{(x-i\varepsilon)(x-i\varepsilon+w)}
\frac{K_\nu(c(-x+i\varepsilon))}{K_\nu(-x+i\varepsilon)}dx.
\endalign
$$
Then, using (3.6) and writing the right hand side by $\xi_\nu$,
$$
\xi_\nu(z)=e^{-i\pi\nu}K_\nu(z)-i\pi I_\nu(z),
$$
we get
$$
\varPi_1=
\frac w{2\pi i}\int_0^{R\cos\theta_{R,\varepsilon}}
\frac{e^{-(c-1)(x-i\varepsilon)}}{(x-i\varepsilon)(x-i\varepsilon+w)}
\frac{\xi_\nu(c(x-i\varepsilon))}{\xi_\nu(x-i\varepsilon)}dx.
$$
Hence, letting $\gamma_1^0$ be the line in $D$ defined by
$$
\gamma_1^0\,:\,
z=x-i\varepsilon,\,\,0\leqq x\leqq R\cos\theta_{R,\varepsilon},
$$
it holds that
$$
\varPi_1=\frac w{2\pi i}\int_{\gamma_1^0}
\frac{e^{-(c-1) z}}{z(z+w)}
\frac{\xi_\nu(cz)}{\xi_\nu(z)}dz.
$$
Here we define three paths as follows:
$$
\align
&\gamma_1^1\,:\,z=\varepsilon e^{i\theta},
\,\,-\pi/2\leqq\theta\leqq0,\\
&\gamma_1^2\,:\,z=x,\,\,\varepsilon\leqq x\leqq R,\\
&\gamma_1^3\,:\,z=Re^{i\theta},
\,\,-\theta_{R,\varepsilon}\leqq\theta\leqq0.
\endalign
$$
Since $w$ is inside $\gamma$, we have that $|\IM w|>\varepsilon$
if $\RE w<0$.
Recall that there is no zero of $K_\nu$
on the real axis. Then we may apply the Cauchy integral theorem
for the integral on the contour consisting of $\gamma_1^0$,
$\gamma_1^1$, $\gamma_1^2$and $\gamma_1^3$ to obtain
$$
\varPi_1=\varPi_1^1+\varPi_1^2-\varPi_1^3,
$$
where
$$
\split
&\varPi_1^1=\frac w{2\pi}\int_{-\pi/2}^0
\frac{e^{-(c-1)\varepsilon e^{i\theta}}}{\varepsilon e^{i\theta}+w}
\frac{\xi_\nu(c\varepsilon e^{i\theta})}{\xi_\nu(\varepsilon e^{i\theta})}
d\theta,\\
&\varPi_1^2=\frac w{2\pi i}\int_\varepsilon^R
\frac{e^{-(c-1) x}}{x(x+w)}\frac{\xi_\nu(cx)}{\xi_\nu(x)}dx,\\
&\varPi_1^3=\frac w{2\pi}\int_{-\theta_{R,\varepsilon}}^0
\frac{e^{-(c-1) Re^{i\theta}}}{Re^{i\theta}+w}
\frac{\xi_\nu(cR e^{i\theta})}{\xi_\nu(R e^{i\theta})}d\theta.
\endsplit
$$
$\varPi_1^3$ tends to $0$ as $R\to\infty$. In fact,
noting that $\xi_\nu(xe^{i\theta})=K_\nu(xe^{i(\theta+\pi)})$
holds for $x>0$, we obtain from Lemma 3.1
$$
\frac{\xi_\nu(cR e^{i\theta})}{\xi_\nu(R e^{\theta})}
=\frac1{\sqrt c}e^{-(c-1) Re^{i(\theta+\pi)}}
\frac{1+E_1(cRe^{i(\theta+\pi)})}{1+E_1(Re^{i(\theta+\pi)})}
$$
for $|\theta|<\pi/6$, which yields
$$
\biggl| e^{-(c-1) Re^{i\theta}}
\frac{\xi_\nu(cR e^{i\theta})}{\xi_\nu(R e^{\theta})}\biggr|
\leqq\frac1{\sqrt c}\frac{1+C_1/cR}{1-C_1/R}
\leqq C_3
\tag3.11
$$
for large $R$ and a positive constant $C_3$ which is
independent of $R$ and $\theta$.
Since $0<\theta_{R,\varepsilon}<\pi/6$, we see $\varPi_1^3\to0$ as $R\to\infty$.

Furthermore (3.11) shows that the function
$e^{-(c-1)x}\xi_\nu(cx)/\xi_\nu(x)$ is bounded
on $[\varepsilon,\infty)$ and that $\varPi_1^2$ converges as $R\to\infty$.
Therefore it holds that 
$$
\lim_{R\to\infty}\varPi_1
=\frac w{2\pi}\int_{-\pi/2}^0
\frac{e^{-(c-1)\varepsilon e^{i\theta}}}{\varepsilon e^{i\theta}+w}
\frac{\xi_\nu(c\varepsilon e^{i\theta})}
{\xi_\nu(\varepsilon e^{i\theta})}d\theta
+\frac w{2\pi i}\int_\varepsilon^\infty
\frac{e^{-(c-1) x}}{x(x+w)}\frac{\xi_\nu(cx)}{\xi_\nu(x)}dx.
$$

In the same way, we can show that
$$
\lim_{R\to\infty}(-\varPi_3)
=\frac w{2\pi}\int_0^{\pi/2}
\frac{e^{-(c-1)\varepsilon e^{i\theta}}}{\varepsilon e^{i\theta}+w}
\frac{\zeta_\nu(c\varepsilon e^{i\theta})}
{\zeta_\nu(\varepsilon e^{i\theta})}d\theta
-\frac w{2\pi i}\int_\varepsilon^\infty
\frac{e^{-(c-1) x}}{x(x+w)}\frac{\zeta_\nu(cx)}{\zeta_\nu(x)}dx,
$$
where
$$
\zeta_\nu(z)=K_\nu(ze^{-i\pi})=e^{i\pi\nu}K_\nu(z)+i\pi I_\nu(z)
$$
for $z\in D$ (cf. (3.7)). Note that
$$
\frac1{2\pi i}\biggl\{
\frac{\zeta_\nu(cx)}{\zeta_\nu(x)}-\frac{\xi_\nu(cx)}{\xi_\nu(x)}
\biggr\}=
\frac{\cos(\pi\nu)\{I_\nu(cx)K_\nu(x)-I_\nu(x)K_\nu(cx)\}}
{\{K_\nu(x)\}^2+\pi^2\{I_\nu(x)\}^2+2\pi\sin(\pi\nu)K_\nu(x)I_\nu(x)}
$$
and recall that the right hand side is $L_{\nu,c}(x)$.
Then we get
$$
\split
\lim_{R\to\infty}
\varPi(R,\varepsilon)=
&\frac w{2\pi}\int_{-\pi/2}^0
\frac{e^{-(c-1)\varepsilon e^{i\theta}}}{\varepsilon e^{i\theta}+w}
\frac{\xi_\nu(c\varepsilon e^{i\theta})}
{\xi_\nu(\varepsilon e^{i\theta})}d\theta\\
&+\frac w{2\pi}\int_0^{\pi/2}
\frac{e^{-(c-1)\varepsilon e^{i\theta}}}{\varepsilon e^{i\theta}+w}
\frac{\zeta_\nu(c\varepsilon e^{i\theta})}
{\zeta_\nu(\varepsilon e^{i\theta})}d\theta\\
&-\frac w{2\pi}\int_{-\pi/2}^{\pi/2}
\frac{e^{(c-1)\varepsilon e^{i\theta}}}{\varepsilon e^{i\theta}-w}
\frac{K_\nu(c\varepsilon e^{i\theta})}{K_\nu(\varepsilon e^{i\theta})}
d\theta\\
&-\int_\varepsilon^\infty
\frac{we^{-(c-1) x}L_{\nu,c}(x)}{x(x+w)}dx.
\endsplit
\tag3.12
$$

We will calculate the limit of each term of (3.12)
as $\varepsilon\downarrow0$.

\medskip

\proclaim{Lemma 3.2}
Let $c>0$, $\nu\geqq0$ and $|\theta|<\pi$. We have that
$$
\lim_{\varepsilon\downarrow0}
\frac{K_\nu(c\varepsilon e^{i\theta})}{K_\nu(\varepsilon e^{i\theta})}
=\lim_{\varepsilon\downarrow0}
\frac{\xi_\nu(c\varepsilon e^{i\theta})}{\xi_\nu(\varepsilon e^{i\theta})}
=\lim_{\varepsilon\downarrow0}
\frac{\zeta_\nu(c\varepsilon e^{i\theta})}{\zeta_\nu(\varepsilon e^{i\theta})}
=\frac1{c^\nu}.
$$
\endproclaim
\demo{Proof}It is known that
$$
K_\nu(z)=
\cases
\log \biggl(\dfrac2z\biggr)\{1+o(1)\}&\text{if $\nu=0$,}\\
\dfrac{\varGamma(\nu)}2\biggl(\dfrac2z\biggr)^\nu\{1+o(1)\}
\quad&\text{if $\nu>0$}
\endcases
\tag3.13
$$
as $|z|\to0$ in $D$. See \cite{15}, p.111 and \cite{18}, p.512.
Then, it follows from (3.13) that
$$
\frac{K_\nu(c\varepsilon e^{i\theta})}{K_\nu(\varepsilon e^{i\theta})}=
\cases
\dfrac{\log(2/c\varepsilon)-i\theta}{\log(2/\varepsilon)-i\theta}
\dfrac{1+o(1)}{1+o(1)}\quad&\text{if $\nu=0$,}\\
\dfrac1{c^\nu}\dfrac{1+o(1)}{1+o(1)}
\quad&\text{if $\nu>0$,}
\endcases
$$
which converges to $1/c^\nu$ as $\varepsilon\downarrow0$.
Recall the formula
$$
I_\nu(z)=\sum_{n=0}^\infty\frac{(z/2)^{\nu+2n}}{n!\varGamma(n+\nu+1)}
\tag3.14
$$
for $z\in D$ (cf. \cite{18}, p.77). From (3.13) and (3.14)
we deduce that $I_\nu(xe^{i\theta})$
converges and $K_\nu(xe^{i\theta})$ tends to infinity
as $x\downarrow0$. This yields that
$$
\frac{\xi_\nu(c\varepsilon e^{i\theta})}{\xi_\nu(\varepsilon e^{i\theta})}
=\frac{K_\nu(c\varepsilon e^{i\theta})+O(1)}
{K_\nu(\varepsilon e^{i\theta})+O(1)}
=\frac1{c^\nu}\{1+o(1)\}
$$
as $\varepsilon\downarrow0$.
\pagebreak

We can show
$$
\frac{\zeta_\nu(c\varepsilon e^{i\theta})}
{\zeta_\nu(\varepsilon e^{i\theta})}
=\frac1{c^\nu}\{1+o(1)\}
$$
in the same fashion.
\qed
\enddemo

\medskip

The first three terms of the right hand side of (3.12) can be
calculated easily. Indeed, Lemma 3.2 yields that
the first and the second terms converge to $1/4c^\nu$ and
that the third term converges to $1/2c^\nu$.
By (3.13) and (3.14), we can easily see
$$
\frac{L_{\nu,c}(x)}{\cos(\pi\nu)}=
\cases
\dfrac{\log c}{(\log x)^2}\{1+o(1)\}\quad&\text{if $\nu=0$,}\\
\dfrac{c^\nu(1-c^{-2\nu})x^{2\nu}}{2^{2\nu-1}\varGamma(\nu)\varGamma(\nu+1)}
\{1+o(1)\}\quad&\text{if $\nu>0$}
\endcases
\tag3.15
$$
as $x\downarrow0$, which has been noted in \cite{3}, p.29.
Hence the last term
of the right hand side of (3.12) converges as
$\varepsilon\downarrow0$. Therefore we can conclude
$$
\lim_{\varepsilon\downarrow0}\lim_{R\to\infty}
\varPi(R,\varepsilon)=
\frac1{c^\nu}-\int_0^\infty
\frac{we^{-(c-1) x}L_{\nu,c}(x)}{x(x+w)}dx.
$$
Since $K_\mu=K_{-\mu}$ for $\mu\geqq0$, we have that
$K_{\nu+1}(z)=K_{|\nu|+1}(z)$ if $z$ is a zero of $K_\nu$.
Moreover, we can regard $z_{-\nu,j}$
as $z_{\nu,j}$ for $1\leqq j\leqq N(\nu)$.
Therefore we have proven the following.

\medskip

\proclaim{Theorem 3.3}
Let $c>1$, $\nu\in\Bbb R$ and $w$ be a point in $D$ with
$K_\nu(w)\neq0$.

\noindent
\text{\rm{(1)}} If $\nu=\pm1/2$, we have
$$
\frac{K_\nu(cw)}{K_\nu(w)}=\frac{e^{-(c-1)w}}{c^{|\nu|}}.
$$
\text{\rm{(2)}} If $|\nu|<3/2$ and $\nu\neq\pm1/2$, we have
$$
\frac{K_\nu(cw)}{K_\nu(w)}=\frac{e^{-(c-1)w}}{c^{|\nu|}}
-e^{-(c-1)w}\int_0^\infty\frac{we^{-(c-1)x}
L_{|\nu|,c}(x)}{x(x+w)}dx.
$$
\text{\rm{(3)}} If $\nu-1/2$ is an integer and $\nu\neq\pm1/2$,
$$
\frac{K_\nu(cw)}{K_\nu(w)}=\frac{e^{-(c-1)w}}{c^{|\nu|}}
-e^{-(c-1)w}\sum_{j=1}^{N(\nu)}
\frac{we^{(c-1)z_{\nu,j}}}{z_{\nu,j}(w-z_{\nu,j})}
\frac{K_\nu(cz_{\nu,j})}{K_{\nu+1}(z_{\nu,j})}.
$$
\text{\rm{(4)}} If $\nu-1/2$ is not an integer and $|\nu|>3/2$,
$$
\split
\frac{K_\nu(cw)}{K_\nu(w)}=\frac{e^{-(c-1)w}}{c^{|\nu|}}
&-e^{-(c-1)w}\sum_{j=1}^{N(\nu)}
\frac{we^{(c-1)z_{\nu,j}}}{z_{\nu,j}(w-z_{\nu,j})}
\frac{K_\nu(cz_{\nu,j})}{K_{\nu+1}(z_{\nu,j})}\\
&-e^{-(c-1)w}\int_0^\infty\frac{we^{-(c-1)x}
L_{|\nu|,c}(x)}{x(x+w)}dx.
\endsplit
$$
\endproclaim

\medskip

We are ready to complete our proof of Theorem 2.2. We have
$$
\L{G_{a,b}^{(\nu)}}(\lambda)=\frac1{\alpha^\nu}
\frac{K_\nu(\alpha\sqrt\lambda)}{\lambda K_\nu(\sqrt\lambda)},
\quad
\alpha=\frac ab>1.
$$
We need to derive the inverse Laplace transforms
of the following functions:
$$
\align
&p_1(\lambda)=\frac1\lambda e^{-(\alpha-1)\sqrt\lambda},\\
&p_2(\lambda;z)=\frac1{\sqrt\lambda(\sqrt\lambda-z)}
e^{-(\alpha-1)\sqrt\lambda},\quad z\in\Bbb C.\\
\endalign
$$
The results may be well known (cf. \cite{16}),
but we deduce them from the formula
$$
\int_0^\infty e^{-\frac{x^2}{u^2}-\frac{v^2}{x^2}}dx
=\frac{u\sqrt\pi}2e^{-\frac{2v}u},
\quad u,v>0.
\tag3.16
$$

At first, put
$$
q_1(t)=\frac1{2\sqrt{\pi t^3}}\int_{\alpha-1}^\infty
\xi\{\xi-(\alpha-1)\}e^{-\frac{\xi^2}{4t}}d\xi
=\frac1{\sqrt{\pi t}}\int_{\alpha-1}^\infty e^{-\frac{\xi^2}{4t}}d\xi.
$$
Then we get by (3.16)
$$
\split
\int_0^\infty e^{-\lambda t}q_1(t)dt
&=\frac1{2\sqrt\pi}\int_{\alpha-1}^\infty
\xi\{\xi-(\alpha-1)\}\biggl[
\int_0^\infty e^{-\lambda t-\frac{\xi^2}{4t}}t^{-\frac32}dt
\biggr]d\xi\\
&=\frac1{\sqrt\pi}\int_{\alpha-1}^\infty
\xi\{\xi-(\alpha-1)\}\biggl[
\int_0^\infty e^{-\frac{\xi^2x^2}4-\frac\lambda{x^2}}dx
\biggr]d\xi\\
&=\int_{\alpha-1}^\infty \{\xi-(\alpha-1)\}e^{-\sqrt\lambda \xi}d\xi\\
&=\frac1\lambda e^{-(\alpha-1)\sqrt\lambda}
\endsplit
$$
and
$$
\IL{p_1}(t)=q_1(t).
$$

Next we put
$$
q_2(t;z)=\frac1{\sqrt{\pi t}}\int_{\alpha-1}^\infty
e^{-\frac{\xi^2}{4t}+z\{\xi-(\alpha-1)\}}d\xi.
$$
Then we obtain from (3.16)
$$
\split
\int_0^\infty e^{-\lambda t}q_2(t;z)dt
&=\frac1{\sqrt\pi}\int_{\alpha-1}^\infty
e^{z\{\xi-(\alpha-1)\}}\biggl[
\int_0^\infty e^{-\lambda t-\frac{\xi^2}{4t}}t^{-\frac12}dt
\biggr]d\xi\\
&=\frac2{\sqrt\pi}\int_{\alpha-1}^\infty
e^{z\{\xi-(\alpha-1)\}}\biggl[
\int_0^\infty e^{-\lambda x^2-\frac{\xi^2}{4x^2}}dx
\biggr]d\xi\\
&=\frac1{\sqrt\lambda}\int_{\alpha-1}^\infty
e^{z\{\xi-(\alpha-1)\}-\sqrt\lambda \xi}d\xi\\
&=\frac1{\sqrt\lambda (\sqrt\lambda-z)} e^{-(\alpha-1)\sqrt\lambda}.
\endsplit
$$
Hence we get
$$
\IL{p_2}(t)=q_2(t;z).
$$

Now we have shown, for example, for the fourth case where
$\nu-1/2$ is not an integer and $|\nu|>3/2$,
$$
\split
G_{a,b}^{(\nu)}(t)=
&\frac1{\alpha^{\nu+|\nu|}}\frac1{\sqrt{\pi t}}
\int_{\alpha-1}^\infty e^{-\frac{\xi^2}{4t}}d\xi\\
&-\frac1{\alpha^\nu}\frac1{\sqrt{\pi t}}\sum_{j=1}^{N(\nu)}
\frac{K_\nu(\alpha z_{\nu,j})}{z_{\nu,j}K_{\nu+1}(z_{\nu,j})}
\int_{\alpha-1}^\infty e^{-\frac{\xi^2}{4t}+z_{\nu,j}\xi}d\xi\\
&-\frac1{\alpha^\nu}\frac1{\sqrt{\pi t}}\int_{\alpha-1}^\infty
e^{-\frac{\xi^2}{4t}}\biggl[
\int_0^\infty \frac{L_{|\nu|,\alpha}(x)}x e^{-x\xi}dx
\biggr]d\xi.
\endsplit
$$
Finally, a simple change of variables from $\xi$ to $s$
given by $\xi=(a-b)\sqrt{2t/s}$ gives us the formula
in Theorem 2.2 (4). The other cases are simpler.


\head
4. The tail probability of the first hitting time.
\endhead

\noindent
As an application of Theorem 2.2, we show the 
asymptotic behavior of $P(\tau_{a,b}^{(\nu)}>t)$
as $t\to\infty$ when $0<b<a$. In Section 2 we showed it
when $\nu=0$ by the Tauberian theorem.
In \cite{3}, it is shown that, if $\nu<0$, 
$$
P(\tau_{a,b}^{(\nu)}>t)=c_\nu t^\nu \{1+o(1)\}
$$
holds for some constant $c_\nu$. It should also be noted that,
in \cite{19}, Yamazato has discussed on the tail probability
in a general framework, and some Bessel processes may be treated.
We give an explicit expression for the constant
$c_\nu$.

To make the statement clear, we define two constants 
when $\nu-1/2$ is an integer. Put
$$
\sigma_1^{(\nu)}=\frac{(a-b)^{2|\nu|}}{2|\nu|}.
$$
Moreover we set $\sigma_2^{(\nu)}=0$ if $\nu=\pm1/2$ and
$$
\sigma_2^{(\nu)}=b^{2|\nu|}(2|\nu|-1)!\sum_{j=1}^{N(\nu)}
\frac{K_\nu(az_{\nu,j}/b)}{z_{\nu,j}^{2|\nu|+1}K_{\nu+1}(z_{\nu,j})}
e^{\frac{z_{\nu,j}(a-b)}b}
\sum_{k=0}^{2|\nu|-1}\frac1{k!}
\biggl\{-\frac{z_{\nu,j}(a-b)}b\biggr\}^k,
$$
if otherwise.

\medskip

\proclaim{Theorem 4.1} Let $0<b<a$.

\noindent
\text{\rm{(1)}} If $\nu=0$,
$$
P(\tau_{a,b}^{(0)}>t)=\frac{2\log(a/b)}{\log t}
+o\biggl(\frac1{\log t}\biggr).
$$
\text{\rm{(2)}} If $\nu>0$ and $\nu-1/2$ is an integer,
$$
\split
P(\tau_{a,b}^{(\nu)}>t)=1-\biggl(\frac ba\biggr)^{2\nu}
+\biggl(\frac ba\biggr)^{2\nu}\sqrt{\frac2\pi}
\frac{(-1/2)^{\nu-1/2}}{(\nu-1/2)!}&
\biggl\{\sigma_1^{(\nu)}+\biggl(\frac ab\biggr)^\nu\sigma_2^{(\nu)}\biggr\}
\frac1{t^\nu}\\
&+O\biggl(\frac1{t^{\nu+1}}\biggr).
\endsplit
$$
\text{\rm{(3)}} If $\nu<0$ and $\nu-1/2$ is an integer,
$$
P(\tau_{a,b}^{(\nu)}>t)=\sqrt{\frac2\pi}
\frac{(-1/2)^{-\nu-1/2}}{(-\nu-1/2)!}
\biggl\{\sigma_1^{(\nu)}+\biggl(\frac ba\biggr)^\nu\sigma_2^{(\nu)}\biggr\}
t^\nu+O(t^{\nu-1}).
$$
\text{\rm{(4)}} If $\nu>0$ and $\nu-1/2$ is not an integer,
$$
P(\tau_{a,b}^{(\nu)}>t)=1-\biggl(\frac ba\biggr)^{2\nu}+
\biggl(\frac{b^3}{2a}\biggr)^\nu\biggl\{\biggl(\dfrac ab\biggr)^\nu
-\biggl(\dfrac ba\biggr)^\nu\biggr\}\frac1{\varGamma(1+\nu)t^\nu}
+o\biggl(\frac1{t^\nu}\biggr).
$$
\text{\rm{(5)}} If $\nu<0$ and $\nu-1/2$ is not an integer,
$$
P(\tau_{a,b}^{(\nu)}>t)=
\biggl(\dfrac2{ab}\biggr)^\nu\biggl\{\biggl(\dfrac ba\biggr)^\nu
-\biggl(\dfrac ab\biggr)^\nu\biggr\}\dfrac{t^\nu}{\varGamma(1-\nu)}
+o(t^\nu).
$$
\endproclaim

\medskip

\noindent
{\bf Remark 4.2.}
It seems that (4) and (5) also hold when $\nu-1/2$ is an integer.  
But we do not pursue the identities.

\medskip

Before proving this theorem, we give two lemmas.
The following one is the immediate consequence of
Lemma 4.3 given in \cite{3}. We let $m(\nu)$ be
the greatest integer which is not larger than $|\nu|-1/2$.

\medskip

\proclaim{Lemma 4.3} We assume that $|\nu|\geqq1/2$.
If $\nu-1/2$ is an integer,
we have that, for any $0\leqq m\leqq m(\nu)-1$
$$
\lim_{t\to\infty}t^{m+1/2}P(\tau_{a,b}^{(\nu)}>t)=0.
\tag4.1
$$
If $\nu-1/2$ is not an integer,
we have (4.1) for any $0\leqq m\leqq m(\nu)$.
\endproclaim

\medskip

For $t>0$, $\nu\neq0$ and $z\in\Bbb C^-$ set
$$
\split
&\varPsi_1(t)=1-\int_0^t\frac{a-b}{\sqrt{2\pi s^3}}
e^{-\frac{(a-b)^2}{2s}}ds,\\
&\varPsi_2(t;z)=\int_0^t\frac{a-b}{\sqrt{2\pi s^3}}
e^{-\frac{(a-b)^2}{2s}+\frac{z(a-b)\sqrt t}{b\sqrt s}}ds,\\
&\varPsi_3(t;\nu)=\int_0^t\frac{a-b}{\sqrt{2\pi s^3}}e^{-\frac{(a-b)^2}{2s}}
\biggl[\int_0^\infty\frac{L_{|\nu|,a/b}(x)}x
e^{-\frac{x(a-b)\sqrt t}{b\sqrt s}}dx\biggr]ds.
\endsplit
$$
Theorem 2.2 implies that $P(\tau_{a,b}^{(\nu)}>t)$ is
represented by a linear combination of $\varPsi_i$'s.
Changing variables from $s$ to $u$ given by $(a-b)/\sqrt s=u$,
we have
$$
\split
&\varPsi_1(t)=1-\sqrt{\frac2\pi}\int_{(a-b)/\sqrt t}^\infty
e^{-\frac{u^2}2}du
=\sqrt{\frac2\pi}\int_0^{(a-b)/\sqrt t}e^{-\frac{u^2}2}du,\\
&\varPsi_2(t;z)=\sqrt{\frac2\pi}\int_{(a-b)/\sqrt t}^\infty
e^{-\frac{u^2}2+\frac{z\sqrt t u}b}du,\\
&\varPsi_3(t;\nu)=\sqrt{\frac2\pi}\cos(\pi\nu)\varPsi_3^0(t;\nu),
\endsplit
$$
where
$$
\varPsi_3^0(t;\nu)=
\int_{(a-b)/\sqrt t}^\infty e^{-\frac{u^2}2}\biggl[
\int_0^\infty\frac{L_{|\nu|,a/b}^0(x)}x
e^{-\frac{x\sqrt t u}b}dx\biggr]du
$$
and $L_{|\nu|,a/b}^0(x)=L_{|\nu|,a/b}(x)/\cos(\pi\nu)$.
It is obvious that $L_{|\nu|,a/b}^0$ is positive on $(0,\infty)$
since $I_\nu$ and $K_\nu$ is increasing and decreasing
on $(0,\infty)$, respectively. 
For an integer $m$ with $0\leqq m\leqq m(\nu)$,
we set
$$
\split
&\beta_1(m)=\frac{(a-b)^{2m+1}}{2m+1},\\
&\beta_2^{(\nu)}(m)=\sum_{j=1}^{N(\nu)}
\frac{K_\nu(az_{\nu,j}/b)}{z_{\nu,j}K_{\nu+1}(z_{\nu,j})}
\beta(m;z_{\nu,j}),\\
&\beta_3^{(\nu)}(m)=(2m)! \, b^{2m+1}\sum_{k=0}^{2m}\frac1{k!}
\biggl(\frac{a-b}b\biggr)^k\int_0^\infty
\frac{L_{|\nu|,a/b}(x)}{x^{2m-k+2}}e^{-\frac{x(a-b)}b}dx.
\endsplit
$$
where
$$
\beta(m;z)=-(2m)!\biggl(\frac bz\biggr)^{2m+1}e^{\frac{z(a-b)}b}
\sum_{k=0}^{2m}\frac1{k!}\biggl\{-\frac{z(a-b)}b\biggr\}^k.
$$
It follows from (3.7) and Lemma 3.1 that
$$
L_{|\nu|,a/b}^0(x)=\pi\sqrt{\frac ba}e^{(a/b-3)x}
\{1+o(1)\}
$$
as $x\to\infty$. Moreover, (3.15) yields that
$L_{|\nu|,a/b}^0(x)/x^{2m-k+2}$ is asymptotically equal to
a constant multiple of $1/x^{2m-k+2-2|\nu|}$ for small $x$.
Since
$$
2m-k+2-2|\nu|\leqq 2m(\nu)+2-2|\nu|<1,
$$
we have that the following improper integral converges:
$$
\int_0^\infty
\frac{L_{|\nu|,a/b}(x)}{x^{2m-k+2}}e^{-\frac{x(a-b)}b}dx
=\cos(\pi\nu)\int_0^\infty
\frac{L_{|\nu|,a/b}^0(x)}{x^{2m-k+2}}e^{-\frac{x(a-b)}b}dx.
$$

\medskip

\proclaim{Lemma 4.4}
If $|\nu|\geqq1/2$,
$$
\align
&\varPsi_1(t)=\sqrt{\frac2\pi}\sum_{m=0}^{m(\nu)}
\frac{(-1/2)^m\beta_1(m)}{m!}\frac1{t^{m+1/2}}
+O\biggl(\frac1{t^{m(\nu)+3/2}}\biggr),\tag4.2\\
&
\varPsi_2(t;z)=\sqrt{\frac2\pi}\sum_{m=0}^{m(\nu)}
\frac{(-1/2)^m\beta(m;z)}{m!}\frac1{t^{m+1/2}}
+O\biggl(\frac1{t^{m(\nu)+3/2}}\biggr).\tag4.3
\endalign
$$
If $0<|\nu|<1/2$,
$$
\varPsi_3(t;\nu)=\biggl(\frac{b^2}2\biggr)^{|\nu|}
\biggl\{\biggl(\frac ab\biggr)^{|\nu|}-\biggl(\frac ba\biggr)^{|\nu|}\biggr\}
\frac1{\varGamma(1+|\nu|)t^{|\nu|}}+o\biggl(\frac1{t^{|\nu|}}\biggr).
\tag4.4
$$
If $|\nu|>1/2$ and $\nu-1/2$ is not an integer,
$$
\split
\varPsi_3(t;\nu)=
&\sqrt{\frac2\pi}\sum_{m=0}^{m(\nu)}
\frac{(-1/2)^m\beta_3^{(\nu)}(m)}{m!}\frac1{t^{m+1/2}}\\
&+\biggl(\frac{b^2}2\biggr)^{|\nu|}
\biggl\{\biggl(\frac ab\biggr)^{|\nu|}-\biggl(\frac ba\biggr)^{|\nu|}\biggr\}
\frac1{\varGamma(1+|\nu|)t^{|\nu|}}\\
&+o\biggl(\frac1{t^{|\nu|}}\biggr).
\endsplit
\tag4.5
$$
\endproclaim
\demo{Proof} For $x\geqq0$ let
$$
P^{(\nu)}(x)=e^{-x^2/2}-\sum_{m=0}^{m(\nu)}
\frac1{m!}\biggl(-\frac{x^2}2\biggr)^m.
$$
Note that
$$
|P^{(\nu)}(x)|\leqq
\frac{x^{2m(\nu)+2}}{2^{m(\nu)+1}\{m(\nu)+1\}!}.
$$
Hence we have
$$
\varPsi_1(t)=\sqrt{\frac2\pi}\sum_{m=0}^{m(\nu)}\frac1{m!}
\biggl(-\frac12\biggr)^m\int_0^{(a-b)/\sqrt t}u^{2m}du+
\sqrt{\frac2\pi}\int_0^{(a-b)/\sqrt t}P^{(\nu)}(u)du,
$$
which implies (4.2). Similarly, by the formula
$$
\int_\beta^\infty x^n e^{-\mu x}dx=e^{-\beta\mu}
\sum_{k=0}^n \frac{n!}{k!}\frac{\beta^k}{\mu^{n-k+1}}
\tag4.6
$$
for $\beta>0$ and $\RE \mu>0$ (cf. \cite{8}, p.340),
we immediately get (4.3).

For $|\nu|>1/2$ and $\nu-1/2$ is not an integer we have
$$
\split
\varPsi_3^0(t;\nu)=&\sum_{m=0}^{m(\nu)}\frac1{m!}
\biggl(-\frac12\biggr)^m \int_{(a-b)/\sqrt t}^\infty
u^{2m}\biggl[\int_0^\infty\frac{L_{|\nu|,a/b}^0(x)}x
e^{-\frac{x\sqrt t u}b}dx\biggr]du\\
&+\int_{(a-b)/\sqrt t}^\infty P^{(\nu)}(u)
\biggl[\int_0^\infty\frac{L_{|\nu|,a/b}^0(x)}x
e^{-\frac{x\sqrt t u}b}dx\biggr]du.
\endsplit
$$
Then the first term of the right hand side is equal to
$$
\sum_{m=0}^{m(\nu)}\frac{(-1/2)^m\beta_3^{(\nu)}(m)}{\cos(\pi\nu)m!}
\frac1{t^{m+1/2}}
$$
since
$$
\int_{(a-b)/\sqrt t}^\infty
u^{2m}\biggl[\int_0^\infty\frac{L_{|\nu|,a/b}^0(x)}x
e^{-\frac{x\sqrt t u}b}dx\biggr]du
=\frac{\beta_3^{(\nu)}(m)}{t^{m+1/2}},
$$
which is obtained by the Fubini theorem and (4.6).
We set
$$
\widetilde\varPsi_3^0(t;\nu)=
\int_{(a-b)/\sqrt t}^\infty P^{(\nu)}(u)
\biggl[\int_0^\infty\frac{L_{|\nu|,a/b}^0(x)}x
e^{-\frac{x\sqrt t u}b}dx\biggr]du.
$$
Changing variables from $x$ to $y$ given by $x\sqrt t u/b=y$,
we have
$$
\widetilde\varPsi_3^0(t;\nu)=
\int_0^\infty 1_{\left[\frac{a-b}{\sqrt t},\infty\right)}(u)P^{(\nu)}(u)
\biggl[\int_0^\infty L_{|\nu|,a/b}^0\biggl(\frac{by}{\sqrt t u}\biggr)
\frac{e^{-y}}y dy\biggr]du,
$$
where $1_A$ is the indicator function of $A$.
To see the convergence of
$t^{|\nu|}\widetilde\varPsi_3^0(t;\nu)$
as $t\to\infty$, we need to dominate
$$
t^{|\nu|}1_{\left[\frac{a-b}{\sqrt t},\infty\right)}(u)
|P^{(\nu)}(u)| L_{|\nu|,a/b}^0\biggl(\frac{by}{\sqrt t u}\biggr)
\frac{e^{-y}}y
\tag4.7
$$
by an integrable function which is independent of $t$.
We have that (4.7) is equal to
$$
b^{2|\nu|}
1_{\left[\frac{a-b}{\sqrt t},\infty\right)}(u)
|P^{(\nu)}(u)|\frac{L_{|\nu|,a/b}^0(by/\sqrt tu)}{(by/\sqrt tu)^{2|\nu|}}
y^{2|\nu|-1}u^{-2|\nu|}e^{-y}.
\tag4.8
$$
Since
$$
\frac{L_{|\nu|,a/b}^0(x)}{x^{2|\nu|}}e^{-(a/b-3)x}
$$
is bounded on $(0,\infty)$, we have that (4.8) is dominated by
a constant multiple of
$$
1_{\left[\frac{a-b}{\sqrt t},\infty\right)}(u)
|P^{(\nu)}(u)| e^{\frac{y(a-3b)}{\sqrt t u}-y}
y^{2|\nu|-1}u^{-2|\nu|}.
\tag4.9
$$
We have that, if $a\leqq3b$, 
$$
e^{\frac{y(a-3b)}{\sqrt t u}-y}\leqq e^{-y}
$$
and that, if $a>3b$,
$$
e^{\frac{y(a-3b)}{\sqrt t u}-y}\leqq e^{\frac{y(a-3b)}{a-b}-y}
=e^{-\frac{2b}{a-b}y}
$$
for $u\geqq (a-b)/\sqrt t$. Let
$$
\kappa=\min\biggl\{1,\frac{2b}{a-b}\biggr\}
$$
and hence (4.9) is bounded by
$$
|P^{(\nu)}(u)|u^{-2|\nu|}y^{2|\nu|-1}e^{-\kappa y}.
\tag4.10
$$
To see that $|P^{(\nu)}(u)|u^{-2|\nu|}$
is integrable on $(0,\infty)$,
we note that
$$
|P^{(\nu)}(x)|\leqq C_4 \min\{1,x^2/2\}x^{2m(\nu)}
$$
for some constant $C_4$. Then we get
$$
|P^{(\nu)}(x)|u^{-2|\nu|}\leqq
\cases
C_4 u^{2m(\nu)-2|\nu|+2}\quad&\text{if $0<u\leqq\sqrt2$,}\\
C_4 u^{2m(\nu)-2|\nu|}\quad&\text{if $u>\sqrt2$.}
\endcases
$$
Since
$$
\split
&2m(\nu)-2|\nu|+2>2\biggl(|\nu|-\frac32\biggr)-2|\nu|+2>-1,\\
&2m(\nu)-2|\nu|<2\biggl(|\nu|-\frac12\biggr)-2|\nu|<-1,
\endsplit
$$
we see that the function given by (4.10) is integrable
on $(0,\infty)\times(0,\infty)$.
Applying the dominated convergence theorem, the Fubini theorem
and (3.15), we have that $t^{|\nu|}\widetilde\varPsi_3^0(t;\nu)$
tends to
$$
\frac{b^{2|\nu|}(a/b)^{|\nu|}\{1-(a/b)^{-2|\nu|}\}}
{2^{2|\nu|-1}\varGamma(|\nu|)\varGamma(|\nu|+1)}\!
\int_0^\infty P^{(\nu)}(u)u^{-2|\nu|}du\!
\int_0^\infty e^{-y}y^{2|\nu|-1}dy
\tag4.11
$$
as $t\to\infty$. Since
$$
\int_0^\infty e^{-y} y^{2|\nu|-1}dy=\varGamma(2|\nu|)
=\frac{2^{2|\nu|-1}}{\sqrt\pi}\varGamma(|\nu|)
\varGamma\biggl(\frac12+|\nu|\biggr)
$$
(cf. \cite{15}, p.3), (4.11) coincides with
$$
\frac{b^{2|\nu|}}{\sqrt\pi}
\biggl\{\biggl(\frac ab\biggr)^{|\nu|}-\biggl(\frac ab\biggr)^{-|\nu|}\biggr\}
\frac1{\varGamma(1+|\nu|)}\varGamma\biggl(\frac12+|\nu|\biggr)
\int_0^\infty P^{(\nu)}(u)u^{-2|\nu|}du.
$$
Changing variables from $u$ to $v$ given by $v=u^2/2$,
we have
$$
\int_0^\infty P^{(\nu)}(u)u^{-2|\nu|}du
=\frac1{2^{|\nu|+1/2}}\int_0^\infty
\frac1{v^{|\nu|+1/2}}
\biggl\{e^{-v}-\sum_{m=0}^{m(\nu)}\frac{(-1)^m}{m!}v^m\biggr\}dv,
$$
which is equal to
$$
\frac1{2^{|\nu|+1/2}}\varGamma\biggl(\frac12-|\nu|\biggr)
$$
(cf. \cite{8}, p.361). The formula
$$
\varGamma\biggl(\frac12+|\nu|\biggr)
\varGamma\biggl(\frac12-|\nu|\biggr)
=\frac\pi{\cos(\pi\nu)}
$$
(cf. \cite{15}, p.3) immediately yields
$$
\widetilde\varPsi_3^0(t;\nu)=
\sqrt{\frac\pi2}\frac1{\cos(\pi\nu)}
\biggl(\frac{b^2}2\biggr)^{|\nu|}
\biggl\{\biggl(\frac ab\biggr)^{|\nu|}-\biggl(\frac ba\biggr)^{|\nu|}\biggr\}
\frac1{\varGamma(1+|\nu|)t^{|\nu|}}
+o\biggl(\frac1{t^{|\nu|}}\biggr)
$$
and hence, we have (4.5).

When $0<|\nu|<1/2$, it is enough to consider $\varPsi_3^0(t;\nu)$ directly.
We can easily deduce (4.4) in the same way
as $\widetilde\varPsi_3^0(t;\nu)$ for $|\nu|>1/2$.
The calculation is left to the reader.
\qed
\enddemo

\medskip

We are now ready to prove Theorem 4.1.
We need only to show Theorem 4.1 in the case of $\nu\neq0$.

For simplicity we set $c=b/a$ and
$$
w_j^{(\nu)}=\frac{K_\nu(az_{\nu,j}/b)}{z_{\nu,j}K_{\nu+1}(z_{\nu,j})}
$$
for $1\leqq j\leqq N(\nu)$.
We first consider the case when $\nu-1/2$ is not an integer
and $\nu<0$.
In this case, if $\nu<-3/2$, Theorem 2.2 gives
$$
P(\tau_{a,b}^{(\nu)}>t)=\varPsi_1(t)+c^\nu
\sum_{j=1}^{N(\nu)}w_j^{(\nu)}\varPsi_2(t;z_{\nu,j})+
c^\nu\varPsi_3(t;\nu).
$$
Note that $m(\nu)+1/2<|\nu|<m(\nu)+3/2$.
It follows from Lemma 4.4 that
$$
\split
P(\tau_{a,b}^{(\nu)}>t)=&
\sqrt{\frac2\pi}\sum_{m=0}^{m(\nu)}
\frac{(-1/2)^m}{m!\,t^{m+1/2}}\{ \beta_1(m)
+c^\nu\beta_2^{(\nu)}(m)
+c^\nu\beta_3^{(\nu)}(m)\}\\
&+c^\nu\biggl(\frac{b^2}2\biggr)^{-\nu}(c^\nu-c^{-\nu})
\frac{t^\nu}{\varGamma(1-\nu)}+o(t^\nu).
\endsplit
$$
By virtue of Lemma 4.3, we have that,
for any $0\leqq m\leqq m(\nu)$
$$
\beta_1(m)+c^\nu \beta_2^{(\nu)}(m)+c^\nu \beta_3^{(\nu)}(m)=0.
\tag4.12
$$
This immediately yields
$$
P(\tau_{a,b}^{(\nu)}>t)=\biggl(\frac2{ab}\biggr)^\nu(c^\nu-c^{-\nu})
\frac{t^\nu}{\varGamma(1-\nu)}+o(t^\nu).
\tag4.13
$$
If $-3/2<\nu<0$, Theorem 2.2 gives
$$
P(\tau_{a,b}^{(\nu)}>t)=\varPsi_1(t)+c^\nu\varPsi_3(t;\nu).
$$
In the case of $-3/2<\nu<-1/2$, we have by Lemma 4.4 that
$$
P(\tau_{a,b}^{(\nu)}>t)=
\sqrt{\frac2{\pi t}}\{ \beta_1(0)+c^\nu\beta_3^{(\nu)}(0)\}
+\biggl(\frac2{ab}\biggr)^\nu(c^\nu-c^{-\nu})
\frac{t^\nu}{\varGamma(1-\nu)}+o(t^\nu).
$$
Lemma 4.3 yields
$$
\beta_1(0)+c^\nu\beta_3^{(\nu)}(0)=0,
$$
and hence we have (4.13). In the case of $-1/2<\nu<0$,
Lemma 4.4 gives (4.13) since $\varPsi_1(t)$ is of order $1/\sqrt t$.
We therefore obtain Theorem 4.1 (5).

We next consider the case when $\nu-1/2$ is not an integer
and $\nu>0$. If $\nu>3/2$, it follows from Theorem 2.2 that
$$
P(\tau_{a,b}^{(\nu)}>t)=1-c^{2\nu}+c^{2\nu}\varPsi_1(t)+c^\nu
\sum_{j=1}^{N(\nu)}w_j^{(\nu)}\varPsi_2(t;z_{\nu,j})+
c^\nu\varPsi_3(t;\nu),
$$
which is equal to
$$
\split
1-c^{2\nu}&+\sqrt{\frac2\pi}\sum_{m=0}^{m(\nu)}
\frac{(-1/2)^mc^{2\nu}}{m!\,t^{m+1/2}}\{\beta_1(m)+c^{-\nu}\beta_2^{(\nu)}(m)
+c^{-\nu}\beta_3^{(\nu)}(m)\}\\
&+c^\nu\biggl(\frac{b^2}2\biggr)^\nu(c^{-\nu}-c^\nu)
\frac1{\varGamma(1+\nu)t^\nu}+O\biggl(\frac1{t^\nu}\biggr).
\endsplit
$$
We have that (4.12) is equivalent to
$$
\beta_1(m)+c^{-\nu}\beta_2^{(-\nu)}(m)+c^{-\nu}\beta_3^{(-\nu)}(m)=0
$$
for $0\leqq m\leqq m(\nu)$.
Since $z_{-\nu,j}$ is regarded as $z_{\nu,j}$ and
$K_{-\nu+1}(z_{-\nu,j})=K_{\nu+1}(z_{\nu,j})$
for $1\leqq j\leqq N(\nu)$, we obtain that
$\beta_2^{(\nu)}(m)=\beta_2^{(-\nu)}(m)$.
Moreover it is trivial that $\beta_3^{(\nu)}(m)=\beta_3^{(-\nu)}(m)$,
and hence we have that, for $0\leqq m\leqq m(\nu)$
$$
\beta_1(m)+c^{-\nu}\beta_2^{(\nu)}(m)+c^{-\nu}\beta_3^{(\nu)}(m)=0.
$$
This immediately yields
$$
P(\tau_{a,b}^{(\nu)}>t)=1-c^{2\nu}+
\biggl(\frac{b^3}{2a}\biggr)^\nu(c^{-\nu}-c^\nu)
\frac1{\varGamma(1+\nu)t^\nu}+O\biggl(\frac1{t^\nu}\biggr).
\tag4.14
$$
In the other cases, (4.14) can be derived in a similar way.
Then we finish to prove Theorem 4.1~(4).

We lastly consider the case when $\nu-1/2$ is an integer.
Similarly to the case when $\nu-1/2$ is not an integer,
we can deduce the asymptotic behavior for $\nu>0$
from that for $\nu<0$. Hence we shall treat
only the case of $\nu<0$. If $\nu\neq-1/2$,
Theorem 2.2 and Lemma 4.4
give
$$
\split
P(\tau_{a,b}^{(\nu)}>t)
&=\varPsi_1(t)+c^\nu\sum_{j=1}^{N(\nu)}w_j^{(\nu)}\varPsi_2(t;z_{\nu,j})\\
&=\sqrt{\frac2\pi}\sum_{m=0}^{m(\nu)}
\frac{(-1/2)^m}{m!\,t^{m+1/2}}\{ \beta_1(m)
+c^\nu\beta_2^{(\nu)}(m)\}
+O\biggl(\frac1{t^{m(\nu)+3/2}}\biggr).
\endsplit
$$
Note that  $m(\nu)=-\nu-1/2\geqq1$. It follows from Lemma 4.3 that
$$
\split
P(\tau_{a,b}^{(\nu)}>t)
&=\sqrt{\frac2\pi}\frac{(-1/2)^{m(\nu)}}{\{m(\nu)\}!\,t^{m(\nu)+1/2}}
\{ \beta_1(m(\nu))+c^\nu\beta_2^{(\nu)}(m(\nu))\}
+O(t^{\nu-1})\\
&=\sqrt{\frac2\pi}
\frac{(-1/2)^{-\nu-1/2}}{(-\nu-1/2)!}
\{\sigma_1^{(\nu)}+c^\nu\sigma_2^{(\nu)}\}
t^\nu+O(t^{\nu-1}).
\endsplit
$$
Since
$$
P(\tau_{a,b}^{(\nu)}>t)=\varPsi_1(t)=\sqrt{\frac2{\pi t}}(a-b)
+O\biggl(\frac1{t^{3/2}}\biggr)
$$
if $\nu=-1/2$, we obtain Theorem 4.1~(3).
We finish to prove Theorem 4.1.

\medskip

\noindent
{\bf Acknowledgment.} The authors would like to thank
Professor Michal Ryznar for valuable comments
on the first draft of this work and the anonymous referee
for useful information concerning the hitting times
of the Bessel processes.


\enddocument